\documentclass[reqno]{amsart}

\usepackage{fullpage,amssymb,dsfont}

\usepackage{tikz}

\usepackage{hyperref}

\numberwithin{equation}{section}

\newtheorem{theorem}{Theorem}[section]
\newtheorem{corollary}[theorem]{Corollary}
\newtheorem{lemma}[theorem]{Lemma}
\newtheorem{proposition}[theorem]{Proposition}
\theoremstyle{remark}
\newtheorem{remark}[theorem]{Remark}

\theoremstyle{definition}

\newcommand{\bp}{\begin{proof}}
\newcommand{\ep}{\end{proof}}

\mathchardef\mhyph="2D

\newcommand{\C}{\mathbb{C}}
\newcommand{\N}{\mathbb{N}}
\newcommand{\Z}{\mathbb{Z}}

\newcommand\eq{=}

\newcommand{\eps}{\varepsilon}

\newcommand{\BB}{\mathcal{B}}
\newcommand{\D}{\mathcal{D}}

\newcommand{\U}{\mathcal{U}}

\newcommand{\un}{\mathds{1}}

\newcommand{\li}{\ell^\infty}

\newcommand\Dhat{\hat{\Delta}}

\newcommand\AUF{A_u(F)}
\newcommand\MB{\mathcal{M}}

\newcommand\Ad{\operatorname{Ad}}
\newcommand\dimm{{\dim_{\mathrm{min}}}}
\newcommand{\End}{\operatorname{End}}
\newcommand\Hom{\operatorname{Hom}}
\newcommand\Irr{\operatorname{Irr}}
\newcommand\Mor{\operatorname{Mor}}
\newcommand\Nat{\operatorname{Nat}}
\newcommand\Rep{\operatorname{Rep}}
\newcommand\Tr{\operatorname{Tr}}
\newcommand\tr{\operatorname{tr}}

\newcommand\qb[2]{\left[\genfrac{}{}{0pt}{}{#1}{#2}\right]_q}

\newcommand{\circt}%
{\mathbin{%
\mathchoice
{\ooalign{$\ocircle$\cr\hidewidth\raise-.15ex\hbox{$\scriptstyle\top\mkern1.7mu$}\cr}}% Woronowicz style tensor product, USUAL SIZE
{\ooalign{$\ocircle$\cr\hidewidth\raise-.15ex\hbox{$\scriptstyle\top\mkern1.7mu$}\cr}}% Woronowicz style tensor product, USUAL SIZE
{\ooalign{$\scriptstyle\ocircle$\cr\hidewidth\raise-.12ex\hbox{$\scriptscriptstyle\top\mkern1mu$}\cr}}% Woronowicz style tensor product, SCRIPT SIZE
{\ooalign{$\scriptstyle\ocircle$\cr\hidewidth\raise-.12ex\hbox{$\scriptscriptstyle\top\mkern1mu$}\cr}}% Woronowicz style tensor product, SCRIPT SIZE
}}

\begin{document}

\title{Martin boundaries of the duals of free unitary quantum groups}

\author[S. Malacarne]{Sara Malacarne}

\email{saramal@math.uio.no}

\author[S. Neshveyev]{Sergey Neshveyev}

\email{sergeyn@math.uio.no}

\address{Department of Mathematics, University of Oslo, P.O. Box 1053
Blindern, NO-0316 Oslo, Norway}

\date{May 19, 2018; minor corrections April 20, 2019}

\begin{abstract}
Given a free unitary quantum group $G=\AUF$, with $F$ not a unitary $2$-by-$2$ matrix, we show that the Martin boundary of the dual of $G$ with respect to any $G$-$\hat G$-invariant, irreducible, finite range quantum random walk coincides with the topological boundary defined by Vaes and Vander Vennet. This can be thought of as a quantum analogue of the fact that the Martin boundary of a free group coincides with the space of ends of its Cayley tree.
\end{abstract}

\maketitle

\section*{Introduction}

The study of harmonic functions on trees has a long history. In the
early 1960s Dynkin and Malyutov~\cite{DM} considered nearest
neighbor random walks on free groups and obtained an analogue of the
Poisson formula for them by identifying the Martin boundary of such a
group with the space of ends of its Cayley graph. This result was
then generalized by a number of authors, among others by
Cartier~\cite{C}, who considered nearest neighbor, but not
necessarily homogeneous, random walks on trees. Finite range random
walks were subsequently studied by Derriennic~\cite{D} in the
homogeneous case and by Picardello and Woess~\cite{PW} in general.
In both cases the result was the same as before: the Martin boundary
of a tree coincides with its space of ends. This was later
generalized to hyperbolic graphs by Ancona~\cite{A} who considered
finite range random walks on such graphs and showed that (under mild
assumptions) the corresponding Martin boundaries coincide with the
Gromov boundaries. A related result was also obtained by
Kaimanovich~\cite{K} who studied only bounded harmonic functions on
hyperbolic groups but went beyond finite range random walks.

The natural quantum analogues of free groups are the duals of free quantum groups of Van Daele and Wang~\cite{VDW}. An interpretation of these duals as vertices of some quantum trees was proposed by Vergnioux~\cite{V}, and an analogue of the end compactification was defined by Vaes and Vergnioux~\cite{VV} in the free orthogonal case and by Vaes and Vander Vennet in the free unitary case~\cite{VVV}. The duals of free orthogonal quantum groups behave in many respects as the group $\Z$, so the corresponding tree has valency~$2$, but the quantum boundary is less trivial than this analogy might suggest; for example, for the dual of $SU_q(2)$ ($q\ne\pm1$) it is the Podle\'{s} quantum sphere $S^2_q$. An identification of the quantum Poisson and Martin boundaries with the quantum spaces of ends was obtained in \cite{I1,NT1,VV}, see also the earlier paper by Biane~\cite{B1} on the dual of~$SU(2)$, which to a great extent motivated the development of the quantum theory.

The duals of free unitary quantum groups behave more like the free group with two generators. In this case Vaes and Vander Vennet showed~\cite{VVV} that the Poisson boundary is isomorphic to the quantum boundary they defined, more precisely, that the von Neumann algebra of bounded harmonic functions is isomorphic to the von Neumann algebra generated by the C$^*$-algebra of continuous functions on the boundary in the GNS-presentation defined by a particular state. Classically this would mean a computation of the Martin boundary as a measure space rather than as a topological space, where the measure is defined as the hitting distribution of the random walk.
%Equivalently, this would give an analogue of the Poisson formula for all bounded harmonic functions as opposed to all positive harmonic functions.
The goal of the present paper is to show that the Martin boundary of the dual of a free unitary quantum group does coincide with the boundary defined in~\cite{VVV} topologically, that is, we have an equality of C$^*$-algebras.

It should be added that a precise relation between quantum Poisson and Martin boundaries has not been established. It was shown in~\cite{NT1} that for the dual of $SU_q(2)$ the computation of the Martin boundary allows one to easily identify the Poisson boundary as well. This was also recently confirmed in~\cite{J} by a different argument, and the same is true for the duals of free orthogonal quantum groups. But a general result of this sort is not available yet, so we cannot claim that our result automatically covers the computation of the Poisson boundary in~\cite{VVV}.

If one wants to compute the Martin boundary of the dual a free unitary quantum group $G=\AUF$, it is natural to try to follow the strategy used in~\cite{VV} in the free orthogonal case. However, this requires a more precise information on the Green functions of particular classical random walks on trees than what seems to be available. We propose a different approach, which also reduces the computation to a classical problem (and works equally well in the free orthogonal case). The idea is that by considering the spectral subspaces of~$\ell^\infty(\hat G)$ with respect to the adjoint action of $G$, we see that a quantum random walk is determined by a collection of operators on the $\ell^2$-spaces of branches of the classical tree $\Irr(G)$. It turns out that the matrix coefficients of any two such operators corresponding to the common branch of definition are exponentially small perturbations of each other (with respect to the distance to the root of $\Irr(G)$). Our main technical result roughly says that the corresponding classical Green functions all have the same asymptotics at infinity. While this is probably not surprising, we point out that the tree $\Irr(G)$ has exponential growth, so it is not immediately obvious that even an exponentially small perturbation is enough to compensate for this growth.

%The paper is organized as follows. In Section~\ref{s1} we collect the basic definitions and some known results on quantum random walks and free unitary quantum groups. In %Section~\ref{s2} we prove a few auxiliary results on the boundary of $\hat G$ defined in~\cite{VVV}. These results rely heavily on estimates obtained in~\cite{VV,VVV}. In %Section~\ref{s3} we reduce the problem of identification of the Martin boundary of $\hat G$ to two statements about classical Green functions. Finally, in Section~\ref{s4} we prove %those statements, thus completing the argument.

\medskip

\paragraph{\bf Acknowledgement} {\small The research leading to these results
received funding from the European Research Council under the
European Union's Seventh Framework Programme (FP/2007-2013) / ERC
Grant Agreement no. 307663. Part of it was carried out during the authors'
visit to the Texas A\&M University. The authors are grateful to Michael Brannan, Ken Dykema, Guoliang Yu and the
staff of the university for their hospitality.}

\bigskip

\section{Preliminaries}\label{s1}

\subsection{Compact and discrete quantum groups}

In this section we recall the basic notions related to compact and discrete quantum groups, see, e.g.,~\cite{NT} for more details.

A \emph{compact quantum group} $G$ is given by a unital C$^*$-algebra $C(G)$ and a unital $*$-homomorphism $$\Delta\colon C(G)\to C(G)\otimes C(G)$$ such that $(\Delta \otimes \iota) \Delta =
(\iota \otimes \Delta) \Delta$ and the spaces
$
(C(G) \otimes 1) \Delta(C(G))\ \ \text{and}\ \ (1 \otimes C(G))\Delta(C(G))
$
are dense in $C(G) \otimes C(G)$.

A \emph{finite dimensional unitary representation} of $G$ is a pair $(U,H_U)$ consisting of a finite dimensional Hilbert space $H_U$ and a unitary element $U\in B(H_U)\otimes C(G)$ such that $\Delta(U)=U_{12}U_{13}$. The tensor product of two representations $U$ and $V$ is defined by $U\otimes V=U_{13}V_{23}\in B(H_U)\otimes B(H_V)\otimes C(G)$. With this tensor product,
%(The tensor product is denoted by $\times$ in~\cite{NT}, another commonly used notation is $\circt$.)
the finite dimensional unitary representations of $G$ form a C$^*$-tensor category $\Rep G$.

The space $\C[G]$ of matrix coefficients of all finite dimensional unitary representations of $G$ is a Hopf $*$-algebra with comultiplication $\Delta$ and antipode $S$ such that $(\iota\otimes S)(U)=U^*$ for all $U\in\Rep G$. The dual space $\U(G)=\C[G]^*$ is a $*$-algebra with multiplication dual to $\Delta$, so $\omega\eta=(\omega\otimes\eta)\Delta$, and involution $\omega^*=\bar\omega S$, where $\bar\omega(a)=\overline{\omega(a^*)}$. This is by definition the algebra of all (unbounded) functions on the \emph{dual discrete quantum group} $\Gamma=\hat G$.

More concretely, $\U(G)$ can be described as follows. For every $U\in\Rep G$ we have a $*$-representation $\pi_U\colon \U(G)\to B(H_U)$ defined by $\pi_U(\omega)=(\iota\otimes\omega)(U)$. Consider the set $I=\Irr(G)$ of isomorphism classes of irreducible unitary representations of $G$. For every $s\in I$ choose a representative~$U_s$. We write $\pi_s$, $H_s$, etc., instead of $\pi_{U_s}$, $H_{U_s}$. Then the collection of representations $(\pi_s)_{s\in I}$ defines a $*$-isomorphism $\U(G)\cong\prod_{s\in I}B(H_s)$. The subalgebra $\ell^\infty(\Gamma)\subset\U(G)$ corresponding to $\ell^\infty\text{-}\bigoplus_{s\in I}B(H_s)$ under this isomorphism is, by definition, the algebra of bounded functions on $\Gamma$. We are not going to distinguish between the algebras $\ell^\infty(\Gamma)$ and $\ell^\infty\text{-}\bigoplus_{s\in I}B(H_s)$. We will also consider the subalgebras
$$
c_0(\Gamma)=c_0\text{-}\bigoplus_{s\in I}B(H_s)\ \ \text{and}\ \ c_c(\Gamma)=\bigoplus_{s\in I}B(H_s)
$$
of functions vanishing at infinity and of finitely supported functions, respectively.

By duality the product on $\C[G]$ defines a comultiplication
$$
\Dhat\colon\U(G)\to\U(G\times G):=(\C[G]\otimes\C[G])^*\cong\prod_{s,t\in I}B(H_s)\otimes B(H_t).
$$
By restriction it gives a map
$
\Dhat\colon \ell^\infty(\Gamma)\to \ell^\infty(\Gamma)\bar\otimes \ell^\infty(\Gamma).
$

The antipode $S$ on $\C[G]$ defines an antipode $\hat S$ on $\U(G)$ by $\hat S(\omega)=\omega S$. There exists a unique element $\rho\in\U(G)$, called the \emph{Woronowicz character}, such that $\rho$ is positive, invertible, $\hat S^2=\Ad\rho$ and $\Tr(\pi_U(\rho))=\Tr(\pi_U(\rho^{-1}))$ for all $U\in\Rep G$. It also has the property $\Dhat(\rho)=\rho\otimes\rho$ justifying the name character. The number $\Tr(\pi_U(\rho))$ is called the \emph{quantum dimension} of $U$ and denoted by $\dim_q (U)$. We will write $\dim_q(s)$ instead of $\dim_q(U_s)$.

The element $\rho$ allows one to define a canonical duality functor on $\Rep G$. Namely, given
a finite dimensional unitary representation $U$ of $G$, the \emph{conjugate representation} is defined by
$$
\bar U=(j(\pi_U(\rho^{1/2}))\otimes1)(j\otimes\iota)(U^*)(j(\pi_U(\rho^{-1/2}))\otimes1)\in B(\bar H_U)\otimes C(G),
$$
where $\bar H_U$ is the complex conjugate Hilbert space to $H_U$ and $j$ is the canonical $*$-anti-isomorphism $B(H_U)\to B(\bar H_U)$, $j(T)\bar\xi=\overline{T^*\xi}$. Note that the element $U^c=(j\otimes\iota)(U^*)$ already satisfies $(\iota\otimes\Delta)(U^c)=U^c_{12}U^c_{13}$, and the additional conjugation by $j(\pi_U(\rho^{1/2}))\otimes1$ is needed only to make this element unitary.

The object $\bar U$ is conjugate to $U$ in $\Rep G$, meaning that there exist morphisms $R_U\colon\un\to \bar U\otimes U$ and $\bar R_U\colon\un\to U\otimes \bar U$, where~$\un$ is the trivial representation of $G$ on the one-dimensional space $\C$, solving the \emph{conjugate equations}: the compositions
$$
U \xrightarrow{\iota \otimes R_U} U\otimes \bar U \otimes U \xrightarrow{\bar R^*_U \otimes \iota} U \quad \mbox{and} \quad \bar U \xrightarrow{\iota \otimes \bar R_U} \bar U \otimes U \otimes \bar U \xrightarrow{R^*_U \otimes \iota}\bar U
$$
are the identity morphisms. Using the Woronowicz character $\rho$ we can define such morphisms by
\begin{equation*}\label{eq:conj}
R_U(1)=\sum_i\bar\xi_i\otimes\pi_U(\rho^{-1/2})\xi_i,\ \ \bar R_U(1)=\sum_i\pi_U(\rho^{1/2})\xi_i\otimes\bar\xi_i,
\end{equation*}
where $\{\xi_i\}_i$ is an orthonormal basis in $H_U$. The pair $(R_U,\bar R_U)$ is \emph{standard}, meaning that $\|R_U\|=\|\bar R_U\|$ and the product $\|R_U\|\,\|\bar R_U\|$ is minimal among all possible solutions of the conjugate equations for~$U$ and~$\bar U$. Note that $\|R_U\|=\|\bar R_U\|=\dim_q (U)^{1/2}$.

\subsection{Quantum random walks}

We now recall the main concepts of the boundary theory of quantum random walks, see~\cite{I1,NT1}.

For every $s\in I=\Irr(G)$, consider the state $\phi_s$ on $B(H_s)$ defined by
$$
\phi_s=\frac{\Tr(\cdot\,\pi_s(\rho^{-1}))}{\dim_q(s)}.
$$
This is a unique state invariant under the left action of $G$ on $B(H_s)$ given by
\begin{equation*}
B(H_s)\to C(G)\otimes B(H_s),\ \ T\mapsto (U_s)_{21}^*(1\otimes T)(U_s)_{21}.
\end{equation*}
Taken together, these actions define the \emph{left adjoint action} of $G$ on $\ell^\infty(\Gamma)$.

Consider now a probability measure $\mu$ on $I$ and define a normal state $\phi_\mu$ on $\ell^\infty(\Gamma)$ by
$$
\phi_\mu=\sum_{s\in I}\mu(s)\phi_s.
$$
Consider the corresponding convolution Markov operator $P_\mu$ on $\ell^\infty(\Gamma)$:
$$
P_\mu=(\phi_\mu\otimes\iota)\Dhat.
$$
This operator commutes with the left adjoint action of $G$ and the right action by translations of $\Gamma$ on~$\ell^\infty(\Gamma)$, which is defined by $\Dhat$. We therefore say that $P_\mu$ defines a $G$-$\Gamma$-\emph{invariant quantum random walk} on~$\Gamma$.

By the $G$-equivariance, the operator $P_\mu$ preserves the center $\ell^\infty(I)$ of $\ell^\infty(\Gamma)$, hence it defines a classical random walk on the set $I$ with transition probabilities $p_\mu(s,t)$ such that
\begin{equation}\label{eq:pmu}
P_\mu(I_t)I_s=p_\mu(s,t)I_s,
\end{equation}
where $I_s$ denotes the unit of $B(H_s)\subset\ell^\infty(\Gamma)$. Explicitly, the transition probabilities are given by
\begin{equation}\label{eq:pmu2}
p_\mu(s,t)=\sum_{r\in I}\mu(r)m^t_{rs}\frac{\dim_q (t)}{\dim_q(r)\dim_q(s)},
\end{equation}
where $m^t_{rs}$ is the multiplicity of $U_t$ in $U_r\otimes U_s$.
We will only be interested in the case when this classical random walk is \emph{irreducible} and \emph{transient}. The irreducibility assumption means that given $s,t$ there exist $s_0=s,s_1,\dots,s_n=t$ such that $p_\mu(s_0,s_1)\dots p_\mu(s_{n-1},s_n)>0$. Equivalently, every representation $U_t$ is a subrepresentation of $U_{r_1}\otimes\dots\otimes U_{r_n}$ for some $r_1,\dots,r_n\in\operatorname{supp}\mu$. In this case we also say that the measure~$\mu$ and the state $\phi_\mu$ are \emph{generating}. The assumption of transience means that the sum of all such products $p_\mu(s_0,s_1)\dots p_\mu(s_{n-1},s_n)$ (with $s_0=s$ and $s_n=t$) is finite for all $s,t\in I$. In this case we also say that~$\phi_\mu$ and~$\mu$ are transient. Similarly to~\cite{DG} it can be shown that, assuming irreducibility, the random walk is transient if the quantum dimension function is nonamenable, in particular, when $\Gamma$ is either nonamenable or non-Kac (recall that the latter means that $\rho\ne1$). We will prove a related result more useful for our purposes in Lemma~\ref{lem:norm} below.

Assuming irreducibility and transience, the \emph{Green kernel} is defined as the completely positive map
$$
G_\mu\colon c_c(\Gamma)\to \ell^\infty(\Gamma),\ \ G_\mu(a)=\sum^\infty_{n=0}P_\mu^n(a),
$$
and then the \emph{Martin kernel} is defined as the map
$$
K_\mu\colon c_c(\Gamma)\to \ell^\infty(\Gamma),\ \ K_\mu(a)=G_\mu(I_e)^{-1}G_\mu(a),
$$
where $e\in I$ is the point corresponding to the trivial representation $\un$. We remark that the element $G_\mu(I_e)$ is a strictly positive scalar in every matrix block $B(H_s)$, so $G_\mu(I_e)^{-1}$ is a well-defined, but possibly unbounded, function on $I$. The elements $K_\mu(a)=G_\mu(I_e)^{-1}G_\mu(a)$ are, however, bounded. Furthermore, if $a\in c_c(\Gamma)$ is positive and nonzero, then $K_\mu(a)$ is a positive invertible element of $\ell^\infty(\Gamma)$.

We have an involution $s\mapsto \bar s$ on $I$ defined by $\bar U_s\cong U_{\bar s}$. For every probability measure $\mu$ on $I$, we define a new probability measure $\check\mu$ on $I$ by $\check\mu(s)=\mu(\bar s)$. The corresponding transition probabilities satisfy
\begin{equation}\label{eq:dualmu}
p_{\check\mu}(s,t)=\frac{\dim_q(t)^2}{\dim_q(s)^2}\,p_\mu(t,s).
\end{equation}
It follows that $\check\mu$ is generating and transient if and only if $\mu$ has the same properties.

Finally, given a generating transient probability measure $\mu$, the \emph{Martin compactification} of $\Gamma$ is defined as the C$^*$-subalgebra $C(\bar\Gamma_{M,\mu})$ of $\ell^\infty(\Gamma)$ generated by $K_{\check\mu}(c_c(\Gamma))$ and $c_c(\Gamma)$. This C$^*$-algebra is invariant under the left adjoint action of $G$ and the right action of $\Gamma$ by translations. The \emph{Martin boundary} is defined by $C(\partial\Gamma_{M,\mu})=C(\bar\Gamma_{M,\mu})/c_0(\Gamma)$.

\subsection{Module categories and categorical random walks} \label{ss:module}

To every C$^*$-algebra $A$ equipped with a left action of~$G$ one can associate the category of equivariant finitely generated Hilbert modules. We will need this only for unital C$^*$-subalgebras $A\subset\ell^\infty(\Gamma)$, equipped with the action given by the restriction of the adjoint action, in which case this category has the following concrete description~\cite[Section~4]{NY0}.

For $U,V\in\Rep G$, consider the subspace
$$
\D_\Gamma(U,V):=\ell^\infty\text{-}\bigoplus_{s\in I}\Hom_G(H_s\otimes H_U,H_s\otimes H_V)\subset
\ell^\infty\text{-}\bigoplus_{s\in I}B(H_s\otimes H_U,H_s\otimes H_V)=\ell^\infty(\Gamma)\otimes B(H_U,H_V),
$$
and define $\D_A(U,V)$ as the intersection of this space with $A\otimes B(H_U,H_V)$. We can then define a new category $\D_A$ with the same objects as in $\Rep G$, but with the morphism spaces $\D_A(U,V)$. (To be more precise, we also need to complete $\D_A$ with respect to subobjects, but whether we do this or not is not going to play any role in this paper.) The category $\Rep G$ can be considered as a (nonfull) subcategory of $\D_A$, if we identify a morphism $T\colon H_U\to H_V$ in $\Rep G$ with the collection of morphisms $\iota_s\otimes T\colon H_s\otimes H_U\to H_s\otimes H_V$. Furthermore, the functors of tensoring on the right by objects of $\Rep G$ define on $\D_A$ the structure of a right $(\Rep G)$-module category.

The morphism spaces $\D_A(U,V)$ are completely determined by the spaces of the form $\D_A(\un,U)$ thanks to the Frobenius reciprocity isomorphisms
\begin{equation}\label{eq:Frobenius}
\D_A(U,V)\cong\D_A(\un,V\otimes\bar U),\ \ T\mapsto (T\otimes\iota_{\bar U})\bar R_U.
\end{equation}
By decomposing $V\otimes\bar U$ into irreducibles we see that $\D_A$ is completely determined already by the spaces $\D_A(\un,U_s)$, $s\in I$.

The $G$-C$^*$-algebra $A$ can be easily reconstructed from the morphism spaces $\D_A(U,V)$: the elements of the form $(\iota\otimes\omega)(T)$, where $T\in \D_A(U,V)\subset \ell^\infty(\Gamma)\otimes B(H_U,H_V)$ and $\omega\in B(H_U,H_V)^*$, span a dense $*$-subalgebra of $A$. By the previous paragraph it is clear that it suffices to take $T\in\D_A(\un,U_s)$, $s\in I$. For each $s\in I$, the linear span of elements of the form $(\iota\otimes\omega)(T)$, with $T\in\D_A(\un,U_s)$ and $\omega\in B(H_e,H_e\otimes H_s)^*=H_s^*$, is nothing other than the spectral subspace of $A$ corresponding to $\bar U_s$.

\smallskip

Assume now that we are given a probability measure $\mu$ on $I$ and consider the corresponding Markov operator $P_\mu$ on $\ell^\infty(\Gamma)$. Since $P_\mu$ commutes with left adjoint action of $G$, the operator $P_\mu\otimes\iota$ is well-defined on $\D_\Gamma(U,V)\subset\ell^\infty(\Gamma)\otimes B(H_U,H_V)$. We will denote this operator by $P_\mu$ again. This operator can be described in purely categorical terms as follows.

The space $\D_\Gamma(U,V)$ can be identified with the space $\Nat_b(\iota\otimes U,\iota\otimes V)$ of bounded natural transformations between the functors $\iota\otimes U$ and $\iota\otimes V$ on $\Rep G$, that is, uniformly bounded collections $(\eta_W)_W$ of morphisms $\eta_W\colon W\otimes U\to W\otimes V$ that are natural in $W$, since any such collection is completely determined by the morphisms $\eta_s=\eta_{U_s}$. For every $s\in I$, we have a contraction $P_s$ on $\Nat_b(\iota\otimes U,\iota\otimes V)$ defined by
$$
P_s(\eta)_W=(\tr_s\otimes\iota)(\eta_{U_s\otimes W}).
$$
where $\tr_s\otimes\iota$ denotes the \emph{normalized categorical partial trace} $\Mor(U_s\otimes X, U_s\otimes Y)\to\Mor(X,Y)$ defined by
$$
(\tr_s\otimes\iota)(T)=\dim_q(s)^{-1}(R^*_s\otimes\iota_Y)(\iota_{\bar s}\otimes T)(R_s\otimes\iota_X),
$$
and where $(R_s,\bar R_s)$ is a standard solution of the conjugate equations for $U_s$ and $U_{\bar s}$. Then
$$
P_\mu=\sum_{s\in I}\mu(s)P_s.
$$
This leads to the following categorical description of the Martin compactification $C(\bar\Gamma_{M,\mu})$ in~\cite{J}, see also~\cite{NY} for the initial motivation and~\cite{DRVV} for a precursor of this picture.

Consider the subspace $\Nat_c(\iota\otimes U,\iota\otimes V)\subset \Nat_b(\iota\otimes U,\iota\otimes V)$ consisting of all natural transformations $\eta$ such that $\eta_s=0$ for all but finitely many $s$. In other words, we consider the subspace
$$
\D_{\Gamma,c}(U,V):=\bigoplus_{s\in I}\Hom_G(H_s\otimes H_U,H_s\otimes H_V)\subset \ell^\infty\text{-}\bigoplus_{s\in I}\Hom_G(H_s\otimes H_U,H_s\otimes H_V).
$$
We have well-defined maps
$$
K_\mu\colon \D_{\Gamma,c}(U,V)\to \D_\Gamma(U,V),\ \ K_\mu(\eta)=G_\mu(I_e)^{-1}G_\mu(\eta),
$$
where $G_\mu(\eta)=\sum^\infty_{n=0}P^n_\mu(\eta)$. Then $\D_{C(\bar\Gamma_{M,\mu})}$ is the smallest C$^*$-subcategory of $\D_\Gamma$ containing $\Rep G$, $\D_{\Gamma,c}$ and all morphisms of the form $K_{\check\mu}(\eta)$, where $\eta$ is a morphism in $\D_{\Gamma,c}$. Using again~\eqref{eq:Frobenius} we see that it suffices to take $\eta\in\D_{\Gamma,c}(\un,U_s)$, $s\in I$.

Note that the previous paragraph does not really add any new information on $C(\bar\Gamma_{M,\mu})$. What we are basically saying is that in computing $K_{\check\mu}\colon c_c(\Gamma)\to\ell^\infty(\Gamma)$ it suffices to consider elements $a$ lying in spectral subspaces of $c_c(\Gamma)$ (with respect to the left adjoint action of $G$) corresponding to $U_s$, and then $K_{\check\mu}(a)$ is again an element of the spectral subspace of $\ell^\infty(\Gamma)$ corresponding to $U_s$. But this is clear, since $c_c(\Gamma)$ is the direct sum of its spectral subspaces and the operator $P_\mu$ is $G$-equivariant. Nevertheless, working with the spaces $\D_{\Gamma}(\un,U_s)$ is more convenient than with the spectral subspaces, and the interpretation of elements of $\D_{\Gamma}(\un,U_s)$ as natural transformations between certain functors makes some computations more transparent.

\subsection{Free unitary quantum groups} \label{ss:free}

In this section we recall some properties of free unitary quantum groups introduced in~\cite{VDW} and studied in detail in~\cite{B}.

Fix a natural number $n\ge 2$ and a matrix $F\in\operatorname{GL}_n(\C)$ such that $\Tr(F^*F)=\Tr((F^*F)^{-1})$. The compact \emph{free unitary quantum group} $A_u(F)$ is defined as the universal unital C$^*$-algebra with generators $u_{ij}$, $1\le i,j\le n$, such that the matrices $U=(u_{ij})_{i,j}$ and $FU^cF^{-1}$ are unitary, where $U^c=(u_{ij}^*)_{i,j}$, equipped with the comultiplication
$$
\Delta(u_{ij})=\sum^n_{k=1}u_{ik}\otimes u_{kj}.
$$
Therefore $U$ defines a unitary representation of the quantum group $A_u(F)$ on the $n$-dimensional Hilbert space, called the \emph{fundamental representation}.

The set $I$ of isomorphism classes of irreducible representations of $A_u(F)$ is the free monoid on letters $\alpha$ and $\beta$, with $\alpha$ corresponding to $U$. The involution $s\mapsto\bar s$ is the anti-automorphism of the monoid defined by $\bar\alpha=\beta$ and $\bar\beta=\alpha$.

From now on we will use the conventions of~\cite{VVV} in that we write $x$ instead of $U_x$ for $x\in I$ whenever convenient. The fusion rules for the representations of $A_u(F)$ are given by
$$
x\otimes y\cong\bigoplus_{z\in I:x=x_0z,y=\bar zy_0}x_0y_0.
$$
Therefore if the last letter of $x$ is the same as the first letter of $y$, then $U_x\otimes U_y$ is irreducible and isomorphic to $U_{xy}$. In this case we write $xy=x\otimes y$.

The Woronowicz character $\rho$ of $A_u(F)$ is determined by the property
$$
\pi_U(\rho)=(F^*F)^t\qquad (\text{the transpose of}\ \ F^*F).
$$
Clearly, $\dim_q(\alpha)=\Tr(F^*F)\ge n$, and the equality holds if and only if $F$ is unitary. Let $q\in(0,1]$ be such that
\begin{equation} \label{eq:q}
\Tr(F^*F)=q+q^{-1}.
\end{equation}
Then $q=1$ if and only if $F$ is a unitary $2$-by-$2$ matrix.

An element $x\in I$ is said to be \emph{indecomposable} if there are only trivial decompositions $x=y\otimes z$, that is, we must have $y=e$ or $z=e$. Equivalently, $x$ is an alternating product of $\alpha$ and $\beta$. Every word $x\in I$ can be written as $x_1\otimes\dots\otimes x_n$, where $x_1,\dots,x_n$ are indecomposable. Then
\begin{equation}\label{eq:qdim}
\dim_q(x)=\prod^n_{i=1}[|x_i|+1]_q,
\end{equation}
where $|x_i|$ denotes the length of $x_i$ and $\displaystyle [n]_q=\frac{q^n-q^{-n}}{q-q^{-1}}$.

\section{Topological boundary of the dual of \texorpdfstring{$\AUF$}{Au(F)}}\label{s2}

Consider a free unitary quantum group $G=\AUF$, with $F$ not a unitary $2$-by-$2$ matrix. Recall that this assumption means that the number $q\in(0,1]$ defined by~\eqref{eq:q} is strictly less than $1$. Denote by $\Gamma$ the dual discrete quantum group.

As we discussed in Section~\ref{ss:free}, the set $I$ of isomorphism classes of irreducible representations of $G$ is a free monoid on letters $\alpha$ and $\beta$. The empty word is denoted by $e$.
Consider the tree with vertex set $I$ such that different elements $x$ and $y$ of $I$ are connected by an edge if and only if one of them is obtained from the other by adding (or removing) one letter on the left. Denote by $\bar I$ the corresponding end compactification of~$I$. The elements of~$\bar I$ are words in $\alpha$ and $\beta$ that are either finite or infinite on the left, and the boundary $\partial I=\bar I\setminus I$ is the set of infinite words. The algebra $C(\bar I)$ of continuous functions on $\bar I$ can be identified with the algebra of functions $f\in\ell^\infty(I)$ such that
$$
|f(yx)-f(x)|\to0\ \text{as}\ x\to\infty,\ \text{uniformly in}\ y\in I.
$$

In \cite{VVV}, Vaes and Vander Vennet extended this construction to $\ell^\infty(\Gamma)$ as follows.
(More precisely, they consider words infinite on the right, while in order to be consistent with our conventions, we consider words infinite on the left.)

For all $x,y\in I$, fix an isometry $V(xy,x\otimes y)\in\Mor(xy,x\otimes y)$. Define ucp maps
$$
\psi_{yx,x}\colon B(H_x)\to B(H_{yx})\ \ T\mapsto V(yx,y\otimes x)^*(1\otimes T)V(yx,y\otimes x).
$$
They do not depend on any choices. Define
\begin{equation}\label{eq:VVVboundary}
\BB=\{a\in\ell^\infty(\Gamma) : \|a_{yx}-\psi_{yx,x}(a_x)\|\to0\ \text{as}\ x\to\infty,\ \text{uniformly in}\ y\in I\}.
\end{equation}
By \cite[Theorem~3.2]{VVV}, this is a unital C$^*$-subalgebra of $\ell^\infty(\Gamma)$ containing $c_0(\Gamma)$. It can therefore be considered as the algebra of continuous functions on a compactification of $\Gamma$. The (noncommutative) algebra of continuous functions on the corresponding boundary is defined by $\BB_\infty=\BB/c_0(\Gamma)$. The left  adjoint action of~$G$ and the right action by translations of~$\Gamma$ on $\ell^\infty(\Gamma)$ define actions of $G$ and $\Gamma$ on $\BB$ and $\BB_\infty$.

If we view $\BB$ only as a $G$-C$^*$-algebra, then from our discussion in Section~\ref{ss:module} we immediately get the following.

\begin{proposition} \label{prop:module-cat-boundary}
Consider the module category $\D_\BB$ associated with the left action of $G$ on $\BB$. Then, for all $U,V\in\Rep G$, the morphism space $\D_\BB(U,V)$ consists of all elements
$$
T\in\li\mbox{-}\bigoplus_{x\in I}\Hom_G(H_x\otimes H_U,H_x\otimes H_V)
$$
such that
$
\|T_{yx}- (V(yx,y\otimes x)^*\otimes\iota_V)(\iota_y\otimes T_x)(V(yx,y\otimes x)\otimes\iota_U)\|\to0
$
as $x\to\infty$, uniformly in $y\in I$.
\end{proposition}

In order to understand better the morphism spaces $\D_\BB(U,V)$ we will now make a particular choice of representations $U_x$ and isometries $V(z,x\otimes y)\in\Mor(z,x\otimes y)$ for all $z\prec x\otimes y$, cf.~\cite[Appendix]{VV}.

Fix representatives $U_\alpha$ and $U_\beta$ of $\alpha$ and $\beta$. All other representatives $U_x$ we construct as follows. If $x=x_1\dots x_n\in I$, with $x_i\in\{\alpha,\beta\}$ for all $i$, then, since $\Mor(x,x_1\otimes\dots\otimes x_n)$ is one-dimensional, $H_{x_1}\otimes \dots \otimes H_{x_n}$ contains a unique invariant subspace on which $U_{x_1}\otimes \dots \otimes U_{x_n}$ is irreducible of class $x$. We take this subspace as~$H_x$ and define $U_x$ as the restriction of $U_{x_1}\otimes \dots \otimes U_{x_n}$ to $H_x$. Denote by $p_x$ the projection in $B(H_{x_1}\otimes \dots \otimes H_{x_n})$ with image $H_x$.

Note that the fusion rules and our choice of representatives imply that $H_{xy}$ is a subspace of $H_x\otimes H_y$.
%, that is,
%\begin{equation}\label{eq:ps}
%p_{xy}\le p_x\otimes p_y\ \ \text{for all}\ \ x,y\in I.
%\end{equation}
We take the embedding map $H_{xy}\to H_x\otimes H_y$ as $V(xy,x\otimes y)$. Then $V(xy,x\otimes y)^*$ coincides with the projection $p^{x\otimes y}_{xy}:=p_{xy}|_{H_x\otimes H_y}$.

Next, viewing $U_\beta$ as the dual of $U_\alpha$, fix a standard solution $(R_\alpha,\bar R_\alpha)$ of the conjugate equations for $U_\alpha$. Put $R_\beta=\bar R_\alpha$ and $\bar R_\beta=R_\alpha$. Assume now that $x=x_1\dots x_n$ and $y=y_1\dots y_m$, with $x_i, y_j \in \{\alpha,\beta\}$ for all $i$ and $j$, and assume that $z$ is a subrepresentation of $x\otimes y$, so
$z=x_1\dots x_{n-k}\, y_{k+1} \dots y_m$ for some $k$ and $x_{n-i}=\bar y_{i+1}$ for $i=0,\dots, k-1$. We then define an operator $\tilde V(z,x\otimes y)\colon H_z\to H_x\otimes H_y$ by
$$
\tilde V(z,x\otimes y)=( p_x \otimes p_y)
(\iota ^{\otimes(n-1)} \otimes \bar R_{x_{n-k+1}} \otimes \iota^{\otimes(m-1)})
\dots(\iota^{\otimes(n-k)} \otimes \bar R_{x_n} \otimes \iota^{\otimes(m-k)}).
$$

This can also be described as follows. Our fixed standard solutions $(R_\alpha,\bar R_\alpha)$ and $(R_\beta,\bar R_\beta)=(\bar R_\alpha,R_\alpha)$ of the conjugate equations allow us to construct standard solutions of the conjugate equations for tensor products of $U_\alpha$ and $U_\beta$. By restriction we then get standard solutions $(R_x,\bar R_x)$ of the conjugate equations for all $x\in I$. By construction, we have
$(R_{\bar x},\bar R_{\bar x})=(\bar R_x, R_x)$ and
\begin{equation}\label{eq:stsol}
R_{xy}=(p^{\bar y\otimes\bar x}_{\bar y\bar x}\otimes p^{x\otimes y}_{xy})(\iota_{\bar y}\otimes R_x\otimes\iota_{y})R_y,\ \ \bar R_{xy}=(p^{x\otimes y}_{xy}\otimes p^{\bar y\otimes\bar x}_{\bar y\bar x})(\iota_{x}\otimes \bar R_y\otimes\iota_{\bar x})\bar R_x.
\end{equation}
Then the morphism $\tilde V(z,x\otimes y)$ for $x=sv$, $y=\bar v t$ and $z=st$ is given by
\begin{equation}\label{eq:Vst}
\tilde V(st,sv\otimes \bar v t)=(p^{s\otimes v}_{sv}\otimes p^{\bar v\otimes t}_{\bar v t})(\iota_s\otimes\bar R_v\otimes\iota_t)V(st,s\otimes t).
\end{equation}

Since $\tilde V(st,sv\otimes \bar v t)$ is an element of $\Mor(st,sv\otimes\bar vt)$, it must be isometric up to a scalar factor.

\begin{lemma}\label{lem:Vnorm}
There is a constant $c>0$ depending only on $q$ such that
$$
\dim_q(v)^{1/2}\ge\|\tilde V(st,sv\otimes \bar v t)\|\ge c\dim_q(v)^{1/2}
$$
for all $v,s,t\in I$.
\end{lemma}

\bp
The first inequality is an immediate consequence of~\eqref{eq:Vst}, since $\|\bar R_v\|=\dim_q(v)^{1/2}$.

The second inequality can be deduced from an analogous result for $SU_{-q}(2)$ using arguments similar to those in~\cite{VVV,Ve}.
Namely, we proceed as follows.

We first reduce the proof to a particular case.
If $v=v_1\otimes v_2$, then, using that $\bar R_v= (\iota_{v_1}\otimes \bar R_{v_2}\otimes\iota_{\bar v_1})\bar R_{v_1}$ by~\eqref{eq:stsol}  and that
$
p^{s\otimes v}_{sv}=p^{s\otimes v_1}_{sv_1}\otimes\iota_{v_2}$, $p^{\bar v\otimes t}_{\bar v t}=\iota_{\bar v_2}\otimes  p^{\bar v_1\otimes t}_{\bar v_1 t}$, we see that
$$
\|\tilde V(st,sv\otimes \bar v t)\|=\dim_q(v_2)^{1/2}\|\tilde V(st,sv_1\otimes \bar v_1 t)\|.
$$
We may therefore assume that $v$ is indecomposable. We also assume that $v\ne e$, since $\tilde V(st,s\otimes t)= V(st,s\otimes t)$ is isometric by definition.

Next, write $sv$ as $s_1\otimes s_2v$, with $s_2v$ indecomposable. Assume first that $s_2=e$. Then
$$
\tilde V(st,sv\otimes \bar v t)=(\iota_s\otimes \iota_v\otimes p^{\bar v\otimes t}_{\bar v t})(\iota_s\otimes\bar R_v\otimes\iota_t)V(st,s\otimes t),
$$
hence
$$
\tilde V(st,sv\otimes \bar v t)^*\tilde V(st,sv\otimes \bar v t)=V(st,s\otimes t)^*(\iota_s\otimes(\Tr_{\bar v}\otimes\iota)(p^{\bar v\otimes t}_{\bar v t}))V(st,s\otimes t),
$$
where $\Tr_{\bar v}$ denotes the categorical trace. The morphism $(\Tr_{\bar v}\otimes\iota)(p^{\bar v\otimes t}_{\bar v t})$ is a scalar multiple of $\iota_t$, and applying~$\Tr_t$ we see that this scalar is $\displaystyle\frac{\dim_q(\bar v t)}{\dim_q(t)}$. It follows that
$$
\|\tilde V(st,sv\otimes \bar v t)\|=\frac{\dim_q(\bar v t)^{1/2}}{\dim_q(t)^{1/2}}.
$$
From \eqref{eq:qdim} it is easy to see that this quantity is not smaller than $c\dim_q(v)^{1/2}$ for a constant $c$ depending only on $q$, cf.~\cite[(5)]{VV}.

Assume now that $s_2\ne e$. Then $st=s_1\otimes s_2t$ and $\tilde V(st,sv\otimes \bar v t)=\iota_{s_1}\otimes\tilde V(s_2t,s_2v\otimes \bar v t)$. We see that in this case the computation reduces to the case when $s_1=e$, that is, we may assume that $sv$ is indecomposable. In a similar way we reduce the computation to the case when $\bar vt$ is indecomposable as well, so we assume that all three elements $sv$, $v\bar v$ and $\bar vt$ are indecomposable (and $v\ne e$). In other words, $sv\bar vt$ equals $\alpha\beta\alpha\dots$ or $\beta\alpha\beta\dots$. But in this case the computation of the norm of $\tilde V(st,sv\otimes \bar v t)$ is equivalent to a similar computation for $SU_{-q}(2)$, see~\cite[Lemma~8.7.2]{Ve}, which gives
\begin{equation*}\label{eq:Vnorm}
\|\tilde V(st,sv\otimes \bar v t)\|=\qb{|s|+|t|+|v|+1}{|v|}^{1/2}\qb{|s|+|v|}{|v|}^{-1/2}\qb{|t|+|v|}{|v|}^{-1/2},
\end{equation*}
see~\cite[(7.3)]{VV}, where $\qb{n}{k}$ denote the $q$-binomial coefficients. Using that
$$
\qb{n}{k}=q^{-k(n-k)}\prod^{k-1}_{i=0}\frac{1-q^{2(n-i)}}{1-q^{2(k-i)}}
$$
we easily deduce that $\|\tilde V(st,sv\otimes \bar v t)\|\ge c[|v|+1]_q^{1/2}=c\dim_q(v)^{1/2}$ for a constant $c>0$ depending only on~$q$, cf.~\cite[Lemma~8.5]{VV}.
\ep

In particular, the morphisms $\tilde V(z,x\otimes y)$ are nonzero, so we can define isometries
$$
V(z,x\otimes y)=\frac{\tilde V(z,x\otimes y)}{\|\tilde V(z,x\otimes y)\|}.
$$

The following result will play a crucial role in our computations. Note that it does not depend on the particular choice of isometries $V(z,x\otimes y)$.

\begin{proposition}[{\cite[Lemma~A.1]{VVV}, \cite[Lemma~8.7.3]{Ve}}]\label{prop:VV-estimates}
There is a constant $C$ depending only on $q$ such that
$$
\|(\iota_u\otimes V(z,x\otimes y))p^{u\otimes z}_{uz}-(p^{u\otimes x}_{ux}\otimes\iota_y)(\iota_u\otimes V(z,x\otimes y))\|\le C q^{(|z|+|x|-|y|)/2},
$$
$$
\|(V(z,x\otimes y)\otimes \iota_u)p^{z\otimes u}_{zu}-(\iota_y\otimes p^{y\otimes u}_{yu})(V(z,x\otimes y)\otimes \iota_u)\|\le C q^{(|z|+|y|-|x|)/2},
$$
for all $u,x,y,z\in I$ such that $z\prec x\otimes y$.
\end{proposition}

For our particular choice of $V(z,x\otimes y)$ this implies the following estimates. (In fact, it is not difficult to see that these estimates are equivalent to the ones above, but with different constants.)

\begin{corollary}\label{cor:VV-estimates}
There is a constant $C$ depending only on $q$ such that
$$
\|V(uz,ux\otimes y)p^{u\otimes z}_{uz}-(p^{u\otimes x}_{ux}\otimes\iota_y)(\iota_u\otimes V(z,x\otimes y))\|\le Cq^{(|z|+|x|-|y|)/2},
$$
$$
\|V(zu,x\otimes yu)p^{z\otimes u}_{zu}-(\iota_x\otimes p^{y\otimes u}_{yu})(V(z,x\otimes y)\otimes\iota_u)\|\le Cq^{(|z|+|y|-|x|)/2},
$$
for all $u,x,y,z\in I$ such that $z\prec x\otimes y$.
\end{corollary}

\bp
By construction we have
$$
\tilde V(uz,ux\otimes y)p^{u\otimes z}_{uz}=(p_{ux}^{u\otimes x}  \otimes \iota_y)(\iota_u \otimes \tilde V(z, x \otimes y))p^{u\otimes z}_{uz},
$$
hence
$$
\|\tilde V(z, x \otimes y)\|^{-1}\tilde V(uz,ux\otimes y)p^{u\otimes z}_{uz}=(p_{ux}^{u\otimes x}  \otimes \iota_y)(\iota_u \otimes V(z, x \otimes y))p^{u\otimes z}_{uz},
$$
so in particular
$$
\|\tilde V(z, x \otimes y)\|^{-1}\|\tilde V(uz,ux\otimes y)\|=\|(p_{ux}^{u\otimes x}  \otimes \iota_y)(\iota_u \otimes V(z, x \otimes y))p^{u\otimes z}_{uz}\|.
$$
Therefore in order to prove the first inequality in the formulation it suffices to show that the last norm is close to $1$ up to $C'q^{(|z|+|x|-|y|)/2}$ for some constant $C'$ and that $(p_{ux}^{u\otimes x}  \otimes \iota_y)(\iota_u \otimes V(z, x \otimes y))p^{u\otimes z}_{uz}$ is close to $(p^{u\otimes x}_{ux}\otimes\iota_y)(\iota_u\otimes V(z,x\otimes y))$ up to $C'q^{(|z|+|x|-|y|)/2}$. But this is indeed the case by the previous proposition, since $(\iota_u \otimes V(z, x \otimes y))p^{u\otimes z}_{uz}$ is an isometry on $H_{uz}$.
%since this inequality implies that $\|(p^{u\otimes x}_{ux}\otimes\iota_y)(\iota_u\otimes V(z,x\otimes y))\|\ge 1-C''q^{(|z|+|x|-|y|)/2}$.
The second inequality is proved in a similar way.
\ep

We can now get a description of the morphism spaces $\D_\BB(e,y)$. If $y$ is not of the form $\bar zz$, $z\in I$, then $\Mor(x,x\otimes y)=0$ for all $x$, so $\D_\BB(e,y)=0$. If $y=\bar zz$, then $x\prec x\otimes y$ precisely for $x$ of the form $uz$, $u\in I$. Denote by $\Omega_y$ the set of all such $x$. It forms a branch of the tree $I$ and its closure $\bar\Omega_y$ in $\bar I$ is a clopen subset.

\begin{corollary} \label{cor:Bmod}
Fix $y\in I$ of the form $\bar zz$. Define an element $T=(T_x)_x\in\ell^\infty\text{-}\bigoplus_{x\in I}\Mor(x,x\otimes y)$ by
$$
T_x=\begin{cases}V(x,x\otimes y),& \text{if}\ \ x\in\Omega_y,\\
0,& \text{otherwise}.\end{cases}
$$
Then $\D_\BB(e,y)=C(\bar\Omega_y) T$.
\end{corollary}

\bp
By Proposition~\ref{prop:module-cat-boundary}, in order to show that $T\in \D_\BB(e,y)$ we have to prove that
$$
\|T_{ux}- (p^{u\otimes x}_{ux}\otimes\iota_y)(\iota_u\otimes T_x)V(ux,u\otimes x)\|\to0
$$
as $x\to\infty$, uniformly in $u\in I$. But this is true by the previous corollary. Next, if $S\in \D_\BB(e,y)$, then $S=fT$ for a function $f\in\ell^\infty(\Omega_y)$ such that
$|f(ux)-f(x)|\to 0$ as $x\to\infty$ in $\Omega_y$, uniformly in $u\in I$. But this precisely means that $f\in C(\bar\Omega_y)$.
\ep

Later we will also need the following estimate.

\begin{corollary}\label{cor:VV-estimates2}
There is a constant $C$ depending only on $q$ such that
$$
\|(V(ux,uv\otimes\bar vx)\otimes\iota_y)V(ux,ux\otimes y)-(\iota_{uv}\otimes V(\bar vx,\bar vx\otimes y))V(ux,uv\otimes\bar vx)\|\le Cq^{|x|-|y|/2}
$$
for all $u,v,x,y\in I$ such that $x\prec x\otimes y$.
\end{corollary}

\bp
Using Corollary~\ref{cor:VV-estimates} we see that up to an operator of norm $\le Cq^{|x|-|y|/2}$ the morphism $$(V(ux,uv\otimes\bar vx)\otimes\iota_y)V(ux,ux\otimes y)p^{u\otimes x}_{ux}$$ equals $(V(ux,uv\otimes\bar vx)\otimes\iota_y)(p^{u\otimes x}_{ux}\otimes\iota_y)(\iota_u\otimes V(x,x\otimes y))p^{u\otimes x}_{ux}$. By \eqref{eq:Vst}, the latter operator equals
$$
\|\tilde V(ux,uv\otimes\bar vx)\|^{-1}(p^{u\otimes v}_{uv}\otimes p^{\bar v\otimes x}_{\bar v x}\otimes\iota_y)(\iota_u\otimes\bar R_v\otimes\iota_x\otimes\iota_y)
(p^{u\otimes x}_{ux}\otimes\iota_y)(\iota_u\otimes V(x,x\otimes y))p^{u\otimes x}_{ux}.
$$
By Proposition~\ref{prop:VV-estimates}, this, in turn,
equals
\begin{multline}\label{eq:Vuxy}
\|\tilde V(ux,uv\otimes\bar vx)\|^{-1}(p^{u\otimes v}_{uv}\otimes p^{\bar v\otimes x}_{\bar v x}\otimes\iota_y)(\iota_u\otimes\bar R_v\otimes\iota_x\otimes\iota_y)
(\iota_u\otimes V(x,x\otimes y))p^{u\otimes x}_{ux}\\
= \|\tilde V(ux,uv\otimes\bar vx)\|^{-1}(p^{u\otimes v}_{uv}\otimes p^{\bar v\otimes x}_{\bar v x}\otimes\iota_y)(\iota_u\otimes\bar R_v\otimes V(x,x\otimes y))p^{u\otimes x}_{ux}
\end{multline}
up to an operator of norm not greater than
$$
C\frac{\|\bar R_v\|}{\|\tilde V(ux,uv\otimes\bar vx)\|}q^{|x|-|y|/2}=C\frac{\dim_q(v)^{1/2}}{\|\tilde V(ux,uv\otimes\bar vx)\|}q^{|x|-|y|/2}.
$$
The last quantity is not greater than $Cc^{-1}q^{|x|-|y|/2}$, where $c$ is the constant from Lemma~\ref{lem:Vnorm}.

In a similar way one checks that the right hand side of \eqref{eq:Vuxy} equals $(\iota_{uv}\otimes V(\bar vx,\bar vx\otimes y))V(ux,uv\otimes\bar vx)p^{u\otimes x}_{ux}$ up to an operator of norm not greater than $Cq^{|x|-|y|/2}+Cc^{-1}q^{|x|-|y|/2}$, where $C$ is from Proposition~\ref{prop:VV-estimates} and~$c$~is from Lemma~\ref{lem:Vnorm}.
\ep

\section{Identification of the Martin boundary}\label{s3}

The following is our main result.

\begin{theorem}\label{thm:MartinFreeUnitary}
Consider a free unitary quantum group $G=\AUF$, with $F$ not a unitary $2$-by-$2$ matrix, and a generating finitely supported probability measure $\mu$ on $I=\Irr(G)$. Then the Martin compactification $C(\bar\Gamma_{M,\mu})$ of the discrete quantum group $\Gamma=\hat G$ with respect to the quantum random walk defined by $\mu$ coincides with the compactification $\BB$ defined in~\eqref{eq:VVVboundary}. It follows that the Martin boundary $C(\partial\Gamma_{M,\mu})$ coincides with~$\BB_\infty$.
\end{theorem}

In order to simplify the notation, let us denote the Martin compactification $C(\bar\Gamma_{M,\mu})$ by $\MB$. Since both~$\MB$ and~$\BB$ are $G$-C$^*$-subalgebras of $\ell^\infty(\Gamma)$, with the actions of $G$ coming from the left adjoint action of $G$ on~$\ell^\infty(\Gamma)$, by the discussion in Section~\ref{ss:module} in order to show that they coincide it suffices to check that $\D_\MB(e,y)=\D_\BB(e,y)$ for all $y\in I$.

\smallskip

Consider first $y=e$. Then, by definition, $\D_\BB(e,e)$ coincides with $C(\bar I)$ considered as a subalgebra of the center $\ell^\infty(I)$ of $\ell^\infty(\Gamma)$. On the other hand, $\D_\MB(e,e)$ contains the algebra of continuous functions on the Martin compactification of $I$ defined by the Markov operator with transition probabilities $p_\mu(s,t)$ defined by~\eqref{eq:pmu}. More precisely, denoting by $P$ this Markov operator and by $G_P$ the corresponding classical Green kernel,
$$
G_P(s,t)=\sum^\infty_{n=0}p_\mu^{(n)}(s,t)=\delta_{st}+\sum^\infty_{n=1}\sum_{s=s_0,s_1,\dots,s_n=t}p_\mu(s_0,s_1)\dots p_\mu(s_{n-1},s_n),
$$
in view of \eqref{eq:dualmu} we have
$$
G_{\check\mu}(I_s)I_t=\frac{\dim_q(s)^2}{\dim_q(t)^2}G_P(s,t).
$$
Therefore, denoting by $K_P$ the classical Martin kernel, $\displaystyle K_P(s,t)=\frac{G_P(s,t)}{G_P(e,t)}$, we see that the function $K_P(s,\cdot)\in\ell^\infty(I)$ differs from $K_{\check\mu}(I_s)$ only by a scalar factor:
\begin{equation*}
K_{\check\mu}(I_s)=\dim_q(s)^2 K_P(s,\cdot).
\end{equation*}
Hence the C$^*$-algebra generated by $c_c(I)\subset c_c(\Gamma)$ and the elements $K_{\check\mu}(I_s)$, $s\in I$, coincides with the algebra of continuous functions on the Martin compactification of $I$. (Note that this does not yet exclude the possibility that $\D_\MB(e,e)=\MB\cap\ell^\infty(I)$ is a strictly larger algebra.) By a result of Picardello and Woess~\cite{PW}, see also~\cite[Corollary~26.14]{W} or \cite[Theorem~27.1]{W}, the Martin compactification of $I$ coincides with the end compactification $\bar I$. It is checked in \cite[Section~2]{VVV} that this result can indeed be applied in our case, since the transition probabilities $p_\mu(s,t)$ satisfy the following properties:

\begin{enumerate}
\item[-] the spectral radius $\lim_np^{(n)}_\mu(s,s)^{1/n}$ of the random walk is strictly less than $1$ (see also Remark~\ref{rem:spectralr} below);
\item[-] the random walk has bounded range: there is $S\in\N$ such that $p_\mu(s,t)=0$ whenever $d(s,t)>S$;
\item[-] the random walk is uniformly irreducible (in the sense used in~\cite{W}): there are $\eps_0>0$ and $K\in\N$ such that for any $s,t\in I$ with $d(s,t)=1$ we have $p^{(k)}_\mu(s,t)\ge\eps_0$ for some $k\le K$.
\end{enumerate}
We thus have
\begin{equation}\label{eq:MvsB1}
\D_\BB(e,e)=C(\bar I)=C^*\big(c_c(I), K_{\check\mu}(c_c(I))\big)\subset \D_\MB(e,e).
\end{equation}

\smallskip

Assume next $y=\bar zz$ for some $z\ne e$. Then, by Corollary~\ref{cor:Bmod}, we have $\D_\BB(e,y)=C(\bar\Omega_y)T$, where $T_x=V(x,x\otimes y)$ for $x\in\Omega_y$. For every $t\in\Omega_y$ define $\eta^{(t)}\in\D_{\Gamma,c}(e,y)=\Nat_c(\iota,\iota\otimes y)$ by
\begin{equation}\label{eq:etas}
\eta^{(t)}_s=\begin{cases}V(t,t\otimes y),& \text{if}\ s=t,\\
0,& \text{otherwise}.\end{cases}
\end{equation}
Define numbers $\check q_\mu(s,t)$ for $s,t\in\Omega_y$ by
\begin{equation*}
P_{\check\mu}(\eta^{(t)})_s=\check q_\mu(s,t)V(s,s\otimes y).
\end{equation*}
Motivated by \eqref{eq:dualmu} we then put
$$
q_\mu(s,t)=\frac{\dim_q(t)^2}{\dim_q(s)^2}\,\check q_\mu(t,s).
$$

\begin{lemma} \label{lem:qvsp1}
The numbers $q_\mu(s,t)$ are real and $|q_\mu(s,t)|\le p_\mu(s,t)$ for all $s,t\in\Omega_y$.
\end{lemma}

\bp
Equivalently, we have to show that the numbers $\check q_\mu(s,t)$ are real and $|\check q_\mu(s,t)|\le p_{\check\mu}(s,t)$.

Since $\check\mu$ is a convex combination of point masses, it suffices to consider $\check\mu=\delta_u$. Fix $s,t\in\Omega_y$. We may assume that $t\prec u\otimes s$, as otherwise $\check q_\mu(s,t)=p_{\check\mu}(s,t)=0$. Then, with $\eta^{(t)}$ given by \eqref{eq:etas},
$$
P_{\check\mu}(\eta^{(t)})_s=(\tr_u\otimes\iota)(\eta^{(t)}_{u\otimes s})=(\tr_u\otimes\iota)((V(t,u\otimes s)\otimes\iota_y)V(t,t\otimes y)V(t,u\otimes s)^*),
$$
and hence
\begin{equation}\label{eq:qmu}
\check q_\mu(s,t)=\tr_{u\otimes s}((\iota_u\otimes V(s,s\otimes y)^*)(V(t,u\otimes s)\otimes\iota_y)V(t,t\otimes y)V(t,u\otimes s)^*).
\end{equation}

Somewhat informally, the statement that $\check q_\mu(s,t)$ is real follows now from the fact that in order to compute the above trace we need only to use that $R_\alpha$ and $\bar R_\alpha$ solve the conjugate equations and satisfy $R_\alpha^*R_\alpha=q+q^{-1}=\bar R_\alpha^*\bar R_\alpha$, and all these relations involve only real numbers. This can be formalized, for example, as follows. We can choose orthonormal bases $(e_i)_i$ of $H_\alpha$ and $(f_i)_i$ of $H_\beta$ such that
$$
R_\alpha(1)=\sum_i\lambda_if_i\otimes e_i\ \ \text{and}\ \ \bar R_\alpha(1)=\sum_i\lambda_i^{-1}e_i\otimes f_i
$$
for some $\lambda_i>0$ (with $\sum_i\lambda_i^2=\sum_i\lambda_i^{-2}=q+q^{-1}$). Then the real linear spans of these bases define real forms of $H_\alpha$ and $H_\beta$. Taking tensor products of these forms we get a real form of $H_{x_1}\otimes\dots\otimes H_{x_n}$ for every word $x=x_1\dots x_n$ in $\alpha$ and $\beta$. Since $H_x$ is the orthogonal complement of the images of operators of the form $\iota\otimes\dots\otimes\iota\otimes R_\alpha\otimes\iota\otimes\dots\otimes\iota$ and $\iota\otimes\dots\otimes\iota\otimes \bar R_\alpha\otimes\iota\otimes\dots\otimes\iota$, which respect the real forms of the tensor products, we then also get a real form of $H_x$. Now, in computing the right hand side of \eqref{eq:qmu} we may work only with the real forms of all the spaces involved, hence the result must be real.

\smallskip

In order to prove the second statement of the lemma, note that we also have
$|\check q_\mu(s,t)|=\|P_{\check\mu}(\eta^{(t)})_s\|$. Consider the map $P_{\check\mu}$ on $\Nat_b(\iota\otimes(e\oplus y),\iota\otimes(e\oplus y))\cong\ell^\infty\text{-}\bigoplus_x\End(x\otimes (e\oplus y))$. Then by restricting it to the block corresponding to $t$ and by projecting the image onto the block corresponding to $s$ we get a completely positive map $P_{st}\colon \End(t\otimes (e\oplus y))\to \End(s\otimes (e\oplus y))$. By definition, this map, being divided by~$p_{\check\mu}(s,t)$, becomes unital. Hence $\|P_{st}\|=p_{\check\mu}(s,t)$. Viewing $V(t,t\otimes y)$ as the element
$$
\begin{pmatrix}
0 & 0\\
V(t,t\otimes y) & 0
\end{pmatrix}
\in \End(t\otimes (e\oplus y)),
$$
we conclude that $\|P_{\check\mu}(\eta^{(t)})_s\|\le\|P_{st}\|=p_{\check\mu}(s,t)$.
\ep

Denote by $Q$ the matrix $(q_\mu(s,t))_{s,t\in\Omega_y}$. We can define the ``Green kernel''
$$
G_Q(s,t)=\delta_{st}+\sum^\infty_{n=1}\sum_{s=s_0,s_1,\dots,s_n=t}q_\mu(s_0,s_1)\dots q_\mu(s_{n-1},s_n).
$$
Then
$$
G_{\check\mu}(\eta^{(s)})_t=\frac{\dim_q(s)^2}{\dim_q(t)^2}G_Q(s,t)V(t,t\otimes y).
$$
Therefore if we let  $\displaystyle K_Q(s,t)=\frac{G_Q(s,t)}{G_P(e,t)}$, then
\begin{equation}\label{eq:MartinK}
K_{\check\mu}(\eta^{(s)})=\dim_q(s)^2 K_Q(s,\cdot)T,
\end{equation}
with $T$ as in Corollary~\ref{cor:Bmod}. We thus need to understand the behavior of the functions $K_Q(s,\cdot)$ at infinity. As we will see, this behavior is not that different from that of $K_P(s,\cdot)$. The starting point is the following estimate.

\begin{lemma}\label{lem:qvsp2}
There is a constant $C$ (depending on $q$, $|y|$ and the support of $\mu$) such that
$$
|q_\mu(s,t)-p_\mu(s,t)|\le Cq^{|s|}
$$
for all $s,t\in\Omega_y$.
\end{lemma}

\bp
Recall that $p_\mu(s,t)=0$ if $d(s,t)>S$, and then $q_\mu(s,t)=0$ as well, so we could equally well use~$q^{|t|}$ instead of~$q^{|s|}$ in the formulation. Using \eqref{eq:qdim} it is also easy to see that the ratios $\dim_q(s)/\dim_q(t)$ are bounded when $d(s,t)\le S$. It follows that the statement of the lemma is equivalent to the existence of $C'$ such that $|\check q_\mu(s,t)-p_{\check\mu}(s,t)|\le C'q^{|s|}$.

Similarly to the proof of the previous lemma we may assume that $\check\mu=\delta_u$ and consider $s,t\in\Omega_y$ such that $t\prec u\otimes s$. We may also assume that $|s|-|u|\ge|y|/2$, as there are only finitely many pairs $(s,t)$ as above not satisfying this condition. Then $u=u_0v$, $s=\bar vs_0$ and $t=u_0s_0$ for some $s_0,u_0,v\in I$. Since $|s_0|=|s|-|v|\ge|s|-|u|\ge|y|/2$, we have $s_0\prec s_0\otimes y$. By \eqref{eq:qmu} we have
$$
\check q_\mu(s,t)=\tr_{u\otimes s}((\iota_{u_0v}\otimes V(\bar vs_0,\bar vs_0\otimes y)^*)(V(u_0s_0,u_0v\otimes \bar vs_0)\otimes\iota_y)V(u_0s_0,u_0s_0\otimes y)V(t,u\otimes s)^*).
$$
By Corollary~\ref{cor:VV-estimates2} the last expression is close to
$$
\tr_{u\otimes s}(V(u_0s_0,u_0v\otimes\bar vs_0)V(t,u\otimes s)^*)=\tr_{u\otimes s}(V(t,u\otimes s)V(t,u\otimes s)^*)=\frac{\dim_q(t)}{\dim_q(u)\dim_q(s)}
$$
up to $Cq^{|s_0|-|y|/2}\le Cq^{|s|-|u|-|y|/2}$. By~\eqref{eq:pmu2}, the quantity on the right is exactly $p_{\check\mu}(s,t)$.
\ep

From this we will deduce  in the next section the following result.

\begin{proposition} \label{prop:MartinK1}
For every $s\in\Omega_y$, the function $\displaystyle K_Q(s,\cdot)=\frac{G_Q(s,\cdot)}{G_P(e,\cdot)}$ on $\Omega_y$ extends to a continuous function on $\bar\Omega_y\subset\bar I$.
\end{proposition}

By \eqref{eq:MartinK} this already shows that
\begin{equation}\label{eq:MvsB2}
K_{\check\mu}(\D_{\Gamma,c}(e,y))\subset C(\bar \Omega_y)T=\D_\BB(e,y).
\end{equation}

As opposed to the functions $K_P(s,\cdot)$, it may happen that the functions $K_Q(s,\cdot)$ vanish at some points of the boundary of $\Omega_y$. Nevertheless we have the following result, whose proof we also postpone till the next section.

\begin{proposition}\label{prop:MartinK2}
For every point $t_\infty\in\bar\Omega_y\cap\partial I$, we have $K_Q(s,t_\infty)>0$ for all $s$ sufficiently close to $t_\infty$.
\end{proposition}

It is now easy to complete the proof of the equality $\MB=\BB$.

\bp[Proof of Theorem~\ref{thm:MartinFreeUnitary}]
By \eqref{eq:MvsB1}, \eqref{eq:MvsB2} and our discussion at the end of Section~\ref{ss:module}, we can already conclude that $K_{\check\mu}(c_c(\Gamma))\subset\BB$. Hence $\MB\subset\BB$, since we also have $c_0(\Gamma)\subset\BB$.

Furthermore, $\D_\MB(e,y)$ is a bimodule over $\D_\MB(e,e)$ and by~\eqref{eq:MvsB1} we know that $\D_\MB(e,e)\supset C(\bar I)$. Using Proposition~\ref{prop:MartinK2}, for every function $f\in C(\bar\Omega_y)$ we can find $s_1,\dots,s_n\in\Omega_y$ and $g_1,\dots,g_n\in C(\bar\Omega_y)$ such that $f(t_\infty)=\sum_i g_i(t_\infty)K_Q(s_i,t_\infty)$ for all $t_\infty\in\bar\Omega_y\cap\partial I$. In view of \eqref{eq:MartinK} this shows that, with $T$ as in Corollary~\ref{cor:Bmod}, we have
\begin{equation}\label{eq:MvsB3}
\D_\BB(e,y)=C(\bar\Omega_y)T\subset C(\bar I)K_{\check\mu}(\D_{\Gamma,c}(e,y))+c_0(I)T.
\end{equation}
Hence $\D_\BB(e,y)\subset \D_\MB(e,y)$ and therefore $\BB\subset\MB$.
\ep

\begin{remark}
The equality $\MB=\BB$ and \eqref{eq:MvsB3} imply that the space $C(\bar I)K_{\check\mu}(c_c(\Gamma))+c_0(\Gamma)$ is dense in $\BB=C(\bar\Gamma_{M,\mu})$.
\end{remark}

\section{Estimates of the Green kernels}\label{s4}

In this section we will prove Propositions~\ref{prop:MartinK1} and~\ref{prop:MartinK2}. We thus fix a generating probability measure $\mu$, a point $y\in I$ of the form $\bar zz$, $z\ne e$, and consider the set $\Omega_y\subset I$ of words of the form $uz$. In order to simplify the notation we denote the transition probabilities $p_\mu(s,t)$ by $p(s,t)$.

In the previous section we introduced the numbers $q_\mu(s,t)$, $s,t\in\Omega_y$, which we now denote simply by~$q(s,t)$. The only information we will need about them is that, by Lemmas~\ref{lem:qvsp1} and~\ref{lem:qvsp2}, they are real and
\begin{equation}\label{eq:qvsp}
|q(s,t)|\le p(s,t)\ \ \text{and}\ \ |q(s,t)-p(s,t)|\le Cq^{|s|}\ \ \text{for all}\ \ s,t\in\Omega_y.
\end{equation}
Denote by $P$ and $Q$ the matrices $(p(s,t))_{s,t\in I}$ and $(q(s,t))_{s,t\in\Omega_y}$.

\smallskip

The space $c_c(I)$ can be completed to two different Hilbert spaces, depending on the inner product we choose.
Denote by $\ell^2(I)$ the completion of $c_c(I)$ with respect to the inner product  $(\cdot,\cdot)$
corresponding to the counting measure: $(f,g)=\sum_{x\in I} f(x)\overline{g(x)}$.
Denote by $\ell^2(I,m)$ the completion of $c_c(I)$ with respect to the inner product corresponding to the ``Haar measure'' $m$ given by $m(x)=\dim_q(x)^2$:
$(f,g)_m=\sum_{x\in I}\,m(x)\,f(x)\overline{g(x)}$. For any $f\in c_c(I)$ and $s\in I$, we have
\begin{equation}\label{eq:scalarproducts}
 \left(f,\delta_s\right)=m(s)^{-1}(f,\delta_s)_m.
\end{equation}

\begin{lemma}\label{lem:norm}
The operator $P$ on $c_c(I)$ extends to a bounded operator on $\ell^2(I,m)$ of norm $\lambda<1$.
The operator $Q$ on $c_c(\Omega_y)$ extends to a bounded operator on $\ell^2(\Omega_y,m)$ of norm not greater than $\lambda$.
\end{lemma}

\bp
Denote by $\dimm$ the dimension function on the finite dimensional representations of our quantum group $G=\AUF$ obtained by letting $q=1$  in formula~\eqref{eq:qdim}.

As can be easily seen from~\eqref{eq:pmu2}, the unitary operator $\ell^2(I,m)\to\ell^2(I)$, $\delta_s\mapsto \dim_q(s)\delta_s$, transforms the operator $P$ into the operator
$$
\sum_{r\in I}\frac{\mu(r)}{\dim_q(r)}\lambda_r,
$$
where $\lambda_r$ is the operator of multiplication by $r$ on the fusion algebra $\C[I]$ of $G$ (which we identify with $c_c(I)$ as a space). By~\cite[Lemma~2.7.3]{NT}, the operator $\lambda_r$ on $\ell^2(I)$ is bounded, of norm not greater than $\dimm(r)$. Therefore
$$
\lambda\le \sum_{r\in I}\mu(r)\frac{\dimm(r)}{\dim_q(r)}<1,
$$
since $\dimm(r)<\dim_q(r)$ for all $r\ne e$.

The second statement of the lemma follows from~\eqref{eq:qvsp}, since we get
$$
|(Qf,g)_m|\le (P|f|,|g|)_m
$$
for all $f,g\in c_c(\Omega_y)$.
\ep

\begin{remark}\label{rem:spectralr}
Note that by \eqref{eq:scalarproducts} we have
\begin{equation}\label{eq:pscalarproduct}
p^{(n)}(s,t)=(P^n\delta_t,\delta_s)=m(s)^{-1}(P^n\delta_t,\delta_s)_m.
\end{equation}
From this we see that the spectral radius of the random walk defined by $P$ is not greater than $\lambda$. This (already mentioned) fact that the spectral radius of the random walk is strictly less than one is actually enough for the results below, but since the inequality $\|P\|<1$ slightly simplifies the arguments, we are going to use it.
\end{remark}

For $x\in I$, let us introduce the following subset of $I$:
$$
\Delta_x=\{ux: u\in I\}. %\ \ \Delta^*_x=\{ux':u\in I\}.
$$
We will have to consider the substochastic matrices $P_{\Delta_x}=(p(s,t))_{s,t\in\Delta_x}$ and the corresponding Green kernels, which we denote by $G_{P,\Delta_x}$. For $x\in\Omega_y$, we similarly define $Q_{\Delta_x}$ and $G_{Q,\Delta_x}$.

\smallskip

The Green kernels $G_{P,\Delta_x}$ satisfy the following uniform version of Harnack's inequality.

\begin{lemma}\label{lem:Harnack}
There exists $\delta\in(0,1)$ such that
$$
G_{P,\Delta_x}(s,t)\le \delta^{-d(s,v)}G_{P,\Delta_x}(v,t)\ \ \text{and}\ \ G_{P,\Delta_x}(s,t)\le \delta^{-d(t,v)}G_{P,\Delta_x}(s,v)
$$
for all $x\in I$ and $v,s,t\in\Delta_x$.
\end{lemma}

\bp
As we already mentioned in the previous section, the random walk defined by $P$ is uniformly irreducible. In fact, the following slightly stronger property (also called uniform irreducibility in~\cite{PW})  is shown in~\cite[Section~2]{VVV}. There exist $K\in\N$ and $\delta_0>0$ such that
\begin{enumerate}
\item[-] $p(s,t)\ge\delta_0$ whenever $p(s,t)>0$;
\item[-] for any $s,t\in I$ with $d(s,t)=1$ we have $p^{(k)}(s,t)>0$ for some $k\le K$.
\end{enumerate}
In particular, given $s,t\in I$ with $d(s,t)=1$, there are points $x_0=s,x_1,\dots,x_k=t$, $k\le K$, such that $p(x_i,x_{i+1})\ge\delta_0$ for all $i=0,\dots,k-1$. But then we also have $p(x_ix,x_{i+1}x)\ge\delta_0$ for any $x\in I$. That is, the same constants~$K$ and~$\delta_0$ work for $P_{\Delta_x}$. This implies the lemma; specifically, we can take $\delta=\delta_0^K$.
\ep

The most important property of the Green kernels of random walks on trees (shared also by random walks on hyperbolic graphs, where it is much more difficult to prove~\cite{A}) is their almost multiplicativity along geodesics. We will need the following uniform version of this property, cf.~\cite[Proposition~2.1]{INO}.

For $s,t\in I$, we denote by $[s,t]$ the unique geodesic segment between $s$ and $t$.

\begin{lemma}\label{lem:decompgreenbr}
There exists a constant $C_1>0$ such that
$$
\frac{1}{C_1}\,G_{P,\Delta_x}(s,v)G_{P,\Delta_x}(v,t)\le G_{P,\Delta_x}(s,t)\le C_1\,G_{P,\Delta_x}(s,v)G_{P,\Delta_x}(v,t)
$$
for all $x\in I$, $s,t\in\Delta_x$ and $v\in[s,t]$.
\end{lemma}

\bp
The first inequality is true, with $C_1=(1-\rho)^{-1}$, for any irreducible random walk with spectral radius $\rho<1$. Namely, we have
$$
\frac{G_{P,\Delta_x}(s,v)}{G_{P,\Delta_x}(v,v)}G_{P,\Delta_x}(v,t)\le G_{P,\Delta_x}(s,t),
$$
since the expression on the left is the sum of $p(s_0,s_1)\dots p(s_{n-1},s_n)$ over all paths $(s=s_0,s_1,\dots,s_n=t)$ in~$\Delta_x$ that pass through $v$. As $G_{P,\Delta_x}(v,v)\le G_P(v,v)\le(1-\lambda)^{-1}$ by Remark~\ref{rem:spectralr}, we can take $C_1=(1-\lambda)^{-1}$.

\smallskip

For the second inequality, take $x\in I$, $s,t\in\Delta_x$ and $v\in[s,t]$. Then any path from $s$ to $t$ in $\Delta_x$ contributing to $G_{P,\Delta_x}(s,t)$ must pass through the set $A=\Delta_x\cap B_S(v)$, where $B_S(v)$ is the open ball of radius $S$ with center~$v$. (Recall that $S\in\N$ is such that $p(u,w)=0$ whenever $d(u,w)>S$.) It follows that
$$
G_{P,\Delta_x}(s,t)\le\sum_{u\in A}G_{P,\Delta_x}(s,u)G_{P,\Delta_x}(u,t).
$$
By Lemma~\ref{lem:Harnack}, for every $u\in A$, we have
$$
G_{P,\Delta_x}(s,u)G_{P,\Delta_x}(u,t)\le\delta^{-2(S-1)}G_{P,\Delta_x}(s,v)G_{P,\Delta_x}(v,t).
$$
Since the number of vertices in $B_S(v)$ is not greater than $3\cdot 2^{S-1}-2$, we thus see that for the second inequality in the formulation of the lemma we can take $C_1=3(2/\delta^2)^{S-1}$.
\ep

We are now ready to prove the following key estimate.

\begin{lemma}\label{lem:Gdif}
There is a constant $C_2$ such that
$$
|G_{Q,\Delta_x}(s,t)-G_{P,\Delta_x}(s,t)|\le C_2 q^{|x|}G_{P,\Delta_x}(s,t)
$$
for all $x\in\Omega_y$ and $s,t\in\Delta_x$.
\end{lemma}

\bp
By \eqref{eq:pscalarproduct} we have $G_P(s,t)=m(s)^{-1}((1-P)^{-1}\delta_t,\delta_s)_m$. Similarly,
$$
G_{P,\Delta_x}(s,t)=m(s)^{-1}((1-P_{\Delta_x})^{-1}\delta_t,\delta_s)_m\ \ \text{and}\ \
G_{Q,\Delta_x}(s,t)=m(s)^{-1}((1-Q_{\Delta_x})^{-1}\delta_t,\delta_s)_m.
$$
Define functions $\tilde \delta_u=\delta_u/\sqrt{m(u)}$ for $u\in I$. Then we can write
$$
G_{P,\Delta_x}(s,t)=\frac{\sqrt{m(t)}}{\sqrt{m(s)}}((1-P_{\Delta_x})^{-1}\tilde\delta_t,\tilde\delta_s)_m\ \ \text{and}\ \
G_{Q,\Delta_x}(s,t)=\frac{\sqrt{m(t)}}{\sqrt{m(s)}}((1-Q_{\Delta_x})^{-1}\tilde\delta_t,\tilde\delta_s)_m.
$$
From this we get
$$
|G_{Q,\Delta_x}(s,t)-G_{P,\Delta_x}(s,t)|=\frac{\sqrt{m(t)}}{\sqrt{m(s)}}|((1-Q_{\Delta_x})^{-1}(Q_{\Delta_x}-P_{\Delta_x})(1-P_{\Delta_x})^{-1}\tilde\delta_t,\tilde\delta_s )_m|.
$$
Since $\{\tilde\delta_u\}_{u\in\Delta_x}$ is an orthonormal basis in $\ell^2(\Delta_x,m)$, the scalar product in the expression above equals
$$
\sum_{u,v\in\Delta_x}((1-P_{\Delta_x})^{-1}\tilde\delta_t,\tilde\delta_v)_m((Q_{\Delta_x}-P_{\Delta_x})\tilde\delta_v,\tilde\delta_u)_m((1-Q_{\Delta_x})^{-1}\tilde\delta_u,\tilde\delta_s )_m.
$$
It follows that

\smallskip

$\displaystyle |G_{Q,\Delta_x}(s,t)-G_{P,\Delta_x}(s,t)|$
\begin{align*}
&\le\frac{\sqrt{m(t)}}{\sqrt{m(s)}}\sum_{u,v\in\Delta_x}|((1-P_{\Delta_x})^{-1}\tilde\delta_t,\tilde\delta_v)_m((Q_{\Delta_x}-P_{\Delta_x})\tilde\delta_v,\tilde\delta_u)_m((1-Q_{\Delta_x})^{-1}\tilde\delta_u,\tilde\delta_s )_m|\\
&= \sum_{u,v\in\Delta_x}G_{P,\Delta_x}(v,t)\,|q(u,v)-p(u,v)|\,|G_{Q,\Delta_x}(s,u)|.
\end{align*}
By \eqref{eq:qvsp} we have $|G_{Q,\Delta_x}(s,u)|\le G_{P,\Delta_x}(s,u)$ and
 $|q(u,v)-p(u,v)|\le Cq^{|u|}$. Since we also have $q(u,v)=p(u,v)=0$ if $d(u,v)>S$, we conclude that
$$
|G_{Q,\Delta_x}(s,t)-G_{P,\Delta_x}(s,t)|\le C\sum_{\genfrac{}{}{0pt}{}{u,v\in\Delta_x\colon}{d(u,v)\le S}}q^{|u|}G_{P,\Delta_x}(s,u)G_{P,\Delta_x}(v,t).
$$
Applying Lemma~\ref{lem:Harnack} we then get
$$
|G_{Q,\Delta_x}(s,t)-G_{P,\Delta_x}(s,t)|\le 3\left(\frac{2}{\delta}\right)^{S}C\sum_{u\in\Delta_x}q^{|u|}G_{P,\Delta_x}(s,u)G_{P,\Delta_x}(u,t),
$$
where we used again that every closed ball of radius $S$ contains not more than $3\cdot 2^S-2$ vertices.

In order to estimate the above expression, for every $u\in\Delta_x$ denote by $u'$ the point on the geodesic $[s,t]$ closest to $u$. Then $u'\in[u,s]\cap[u,t]$, so by Lemma~\ref{lem:decompgreenbr} we have
\begin{align}
G_{P,\Delta_x}(s,u)G_{P,\Delta_x}(u,t)&\le C^2_1G_{P,\Delta_x}(s,u')G_{P,\Delta_x}(u',u)
G_{P,\Delta_x}(u,u')G_{P,\Delta_x}(u',t)\nonumber\\
&\le C^3_1 G_{P,\Delta_x}(s,t)G_{P,\Delta_x}(u',u) G_{P,\Delta_x}(u,u'),\label{eq:Hthrice}
\end{align}
and hence
\begin{equation}\label{eq:Gdif}
|G_{Q,\Delta_x}(s,t)-G_{P,\Delta_x}(s,t)|\le 3\left(\frac{2}{\delta}\right)^{S}C\, C_1^3\,G_{P,\Delta_x}(s,t)\sum_{u\in\Delta_x}q^{|u|}G_{P,\Delta_x}(u',u) G_{P,\Delta_x}(u,u').
\end{equation}

It remains to estimate the sum in the above expression on the right. Observe that if $u'\ne x'$ (where $x'\in[s,t]$ is the vertex closest to $x$), then
$$
|u|\ge|u'|=|x'|+d(u',x')\ge|x|+d(u',x').
$$
If $u'=x'$, then the inequality $|u|\ge|x|+d(u',x')$ is obvious. See Figure~\ref{fig} illustrating these two cases.

\begin{figure}
 \begin{center}
 \small
  \begin{tikzpicture}[thin,
   level 1/.style={sibling distance=40mm},
   level 2/.style={sibling distance=25mm},
   level 3/.style={sibling distance=15mm},
   level 4/.style={sibling distance=9mm},
   every circle node/.style={minimum size=1.5mm,inner sep=0mm}]
   \node[circle,fill,label=below:$x$] (root) {} [grow=up]
    child { node [circle,fill] {}
     child { node [circle,fill,label=left:$v$] {}
     }
     child { node [circle,fill,label=left:$v'\eq x'$] {}
      child { node[circle,fill,label=left:$t$] {}
      }
      child { node[circle,fill,label=left:$u'$] {}
       child { node[circle,fill,label=left:$u$] {}
       }
       child {node[circle,fill,label=left:$s$] {}
       }
      }
     }
    };
  \end{tikzpicture}
 \end{center}
 \caption{}\label{fig}
\end{figure}

It follows that
\begin{align}
\sum_{u\in\Delta_x}q^{|u|}G_{P,\Delta_x}(u',u) G_{P,\Delta_x}(u,u')&\le q^{|x|}\sum_{u\in\Delta_x}q^{d(u',x')}G_{P,\Delta_x}(u',u) G_{P,\Delta_x}(u,u')\nonumber\\
&\le q^{|x|}\sum_{u\in\Delta_x}\sum_{w\in[s,t]}q^{d(w,x')}G_{P,\Delta_x}(w,u) G_{P,\Delta_x}(u,w).\label{eq:Gdif1}
\end{align}

Next, for every $w\in[s,t]$, we have
\begin{align}
\sum_{u\in\Delta_x}G_{P,\Delta_x}(w,u) G_{P,\Delta_x}(u,w)&=\sum_{u\in\Delta_w}((1-P_{\Delta_x})^{-1}\tilde\delta_u,\tilde\delta_w)_m((1-P_{\Delta_x})^{-1}\tilde\delta_w,\tilde\delta_u)_m\nonumber\\
&=((1-P_{\Delta_x})^{-2}\tilde\delta_w,\tilde\delta_w)_m\le(1-\lambda)^{-2}.\label{eq:Gdif1.5}
\end{align}
We also have
\begin{equation}\label{eq:Gdif2}
\sum_{w\in[s,t]}q^{d(w,x')}< 2(1+q+q^2+\dots)=2(1-q)^{-1}.
\end{equation}

Combining \eqref{eq:Gdif1}--\eqref{eq:Gdif2}, we get from \eqref{eq:Gdif} that
$$
|G_{Q,\Delta_x}(s,t)-G_{P,\Delta_x}(s,t)|\le 6\left(\frac{2}{\delta}\right)^{S}\frac{C\,C^3_1}{(1-\lambda)^2(1-q)}\,q^{|x|}G_{P,\Delta_x}(s,t),
$$
and we are done.
\ep

\bp[Proof of Proposition~\ref{prop:MartinK1}]
Fix $s\in\Omega_y$ and $t_\infty\in\bar\Omega_y\cap\partial I$. We have to show that as $t\to t_\infty$ the numbers $\displaystyle K_Q(s,t)=\frac{G_Q(s,t)}{G_P(e,t)}$ converge to a finite limit.

Fix $\eps>0$. By Lemma~\ref{lem:Gdif} there exists $x\in\Omega_y$ such that $s\not\in\Delta_x$, $t_\infty\in\bar\Delta_x\cap\partial I$ and
$$
|G_{Q,\Delta_x}(u,t)-G_{P,\Delta_x}(u,t)|\le \eps\, G_{P,\Delta_x}(u,t)\ \ \text{for all}\ \ u,t\in\Delta_x.
$$

Take $t\in\Delta_x$. Then $x\in[s,t]$, so similarly to the proof of Lemma~\ref{lem:decompgreenbr}, every path from $s$ to $t$ in $I$ contributing to $G_P(s,t)$ must pass through the set $A=\Delta_x\cap B_S(x)$ and there is a well-defined moment when it leaves $I\setminus\Delta_x$ and enters $A$ for the last time. It follows that if for $u\in A$ we put
$$
M_P(s,u)=\sum_{n=1}^\infty\sum_{\genfrac{}{}{0pt}{}{s=s_0,s_1,\dots,s_n=u\colon}{s_{n-1\not\in\Delta_x}}}p(s_0,s_1)\dots p(s_{n-1},s_n),
$$
then
\begin{equation}\label{eq:MP}
G_P(s,t)=\sum_{u\in A}M_P(s,u)G_{P,\Delta_x}(u,t).
\end{equation}
We can similarly define $M_Q(s,u)$ and get a decomposition of $G_Q(s,t)$. Then, for any $t,t'\in\Delta_x$, we have
\begin{align*}
|K_Q(s,t)-K_Q(s,t')|&=\left|\sum_{u\in A}\frac{M_Q(s,u)G_{Q,\Delta_x}(u,t)}{G_P(e,t)}-\sum_{u\in A}\frac{M_Q(s,u)G_{Q,\Delta_x}(u,t')}{G_P(e,t')}\right|\\
&\le \eps\sum_{u\in A}\frac{|M_Q(s,u)|G_{P,\Delta_x}(u,t)}{G_P(e,t)}+\eps\sum_{u\in A}\frac{|M_Q(s,u)|G_{P,\Delta_x}(u,t')}{G_P(e,t')}\\
&\qquad\qquad\qquad\qquad +\sum_{u\in A}|M_Q(s,u)|\left|\frac{G_{P,\Delta_x}(u,t)}{G_P(e,t)}-\frac{G_{P,\Delta_x}(u,t')}{G_P(e,t')}\right|\\
&\le\eps(K_P(s,t)+K_P(s,t'))+\sum_{u\in A}|M_Q(s,u)|\left|\frac{G_{P,\Delta_x}(u,t)}{G_P(e,t)}-\frac{G_{P,\Delta_x}(u,t')}{G_P(e,t')}\right|,
\end{align*}
where we used that $|M_Q(s,u)|\le M_P(s,u)$ by~\eqref{eq:qvsp}.
If we could prove that the sum above converges to zero as $t,t'\to t_\infty$, then we would get
$$
\limsup_{t,t'\to t_\infty}|K_Q(s,t)-K_Q(s,t')|\le 2\eps\, K_P(s,t_\infty),
$$
and since $\eps$ could be taken arbitrarily small, we would be able to conclude that the numbers $K_Q(s,t)$ converge to a finite limit as $t\to t_\infty$.

Therefore it remains to show that for every $u\in A$ we have
$$
\lim_{t,t'\to t_\infty}\left|\frac{G_{P,\Delta_x}(u,t)}{G_P(e,t)}-\frac{G_{P,\Delta_x}(u,t')}{G_P(e,t')}\right|=0.
$$
Using again~\eqref{eq:MP}, but now applied to $G_P(e,t)$, we can rewrite this as
$$
\lim_{t,t'\to t_\infty}\left|\left(\sum_{v\in A}M_P(e,v)\frac{G_{P,\Delta_x}(v,t)}{G_{P,\Delta_x}(u,t)}\right)^{-1}-\left(\sum_{v\in A}M_P(e,v)\frac{G_{P,\Delta_x}(v,t')}{G_{P,\Delta_x}(u,t')}\right)^{-1}\right|=0.
$$
In order to prove this, it suffices to show that for all $u,v\in A$ the function $\displaystyle\frac{G_{P,\Delta_x}(v,\cdot)}{G_{P,\Delta_x}(u,\cdot)}$ on $\Delta_x$ extends to a nonvanishing continuous function on $\bar\Delta_x\subset\bar I$. But this is indeed true, since this function is nothing else than the Martin kernel of the random walk on the tree $\Delta_x$ defined by the matrix $P_{\Delta_x}$. Although this matrix is only substochastic, it still satisfies the assumptions of \cite[Corollary~26.14]{W} or \cite[Theorem~27.1]{W}, as we essentially discussed in the proof of Lemma~\ref{lem:Harnack}. So the same arguments as in~\cite{W} show that the Martin compactification of $\Delta_x$ with respect to $P_{\Delta_x}$ coincides with the end compactification, which is~$\bar\Delta_x$.
\ep

\bp[Proof of Proposition~\ref{prop:MartinK2}]
Fix $t_\infty\in\bar\Omega_y\cap\partial I$.
By Lemma~\ref{lem:Gdif} we can find $x\in\Omega_y$ such that $t_\infty\in\bar\Delta_x$ and
$$
G_{Q,\Delta_x}(s,t)\ge\frac{1}{2}G_{P,\Delta_x}(s,t)\ \ \text{for all}\ \ s,t\in\Delta_x.
$$
Suppose next that we can find $u\in\Delta_x$ such that $t_\infty\in\bar\Delta_u$ and
\begin{equation}\label{eq:GDeltaDif}
G_P(s,t)-G_{P,\Delta_x}(s,t)\le\frac{1}{4}G_P(s,t)\ \ \text{for all}\ \ s,t\in\Delta_u.
\end{equation}
Then, by \eqref{eq:qvsp},
$$
|G_Q(s,t)-G_{Q,\Delta_x}(s,t)|\le G_P(s,t)-G_{P,\Delta_x}(s,t)\le\frac{1}{4}G_P(s,t)\le\frac{1}{3}G_{P,\Delta_x}(s,t)\le \frac{2}{3}G_{Q,\Delta_x}(s,t),
$$
whence
$$
G_Q(s,t)\ge\frac{1}{3}G_{Q,\Delta_x}(s,t)\ge\frac{1}{6}G_{P,\Delta_x}(s,t)\ge\frac{1}{8}G_P(s,t),
$$
and therefore
$$
K_Q(s,t_\infty)\ge\frac{1}{8}K_P(s,t_\infty)>0\ \ \text{for all}\ \ s\in\Delta_u.
$$
Therefore it remains to find $u$ such that \eqref{eq:GDeltaDif} is satisfied.

For $s,t\in\Delta_x$, similarly to the proof of Lemma~\ref{lem:decompgreenbr}, the paths from $s$ to $t$ contributing to $G_P(s,t)-G_{P,\Delta_x}(s,t)$ must pass through the set $A=B_S(x)\cap\Delta_x$ of cardinality $2^S-1$. It follows that
$$
G_P(s,t)-G_{P,\Delta_x}(s,t)\le\sum_{v\in A}G_P(s,v)G_P(v,t),
$$
and using  Lemma~\ref{lem:Harnack} we get
$$
G_P(s,t)-G_{P,\Delta_x}(s,t)\le \left(\frac{2}{\delta^2}\right)^SG_P(s,x)G_P(x,t).
$$

Now, take any $u\in\Delta_x$ and assume that $s,t\in\Delta_u$. Let $v$ be the point on $[s,t]$ closest to $u$, or equivalently, to $x$. Then $v\in[s,x]\cap[t,x]$, so similarly to \eqref{eq:Hthrice} we have
$$
G_P(s,x)G_P(x,t)\le C_1^3G_P(s,t)G_P(x,v)G_P(v,x).
$$
Hence
$$
G_P(s,t)-G_{P,\Delta_x}(s,t)\le\left(\frac{2}{\delta^2}\right)^SC_1^3G_P(s,t)G_P(x,v)G_P(v,x).
$$
Therefore for \eqref{eq:GDeltaDif} to be satisfied it suffices to have
\begin{equation}\label{eq:Gproduct}
G_P(x,v)G_P(v,x)\le \frac{1}{4C_1^3}\left(\frac{\delta^2}{2}\right)^S
\end{equation}
for all $v\in\Delta_u$. But, similarly to~\eqref{eq:Gdif1.5},
$$
\sum_{v\in I}G_P(x,v)G_P(v,x)\le(1-\lambda)^{-2}<\infty,
$$
so \eqref{eq:Gproduct} is true for all but finitely many $v$'s in $I$.
\ep

\begin{remark}
Essentially the same proof shows that actually $\displaystyle\frac{K_Q(s,t_\infty)}{K_P(s,t_\infty)}\to1$ as $s\to t_\infty$.
\end{remark}

\bigskip

\end{document}